\RequirePackage{fix-cm}
\documentclass[smallextended,envcountsect]{svjour3}       
\smartqed  
\usepackage{graphicx}
\usepackage{amsmath,amssymb}
\usepackage{array}
\usepackage{enumitem}
\usepackage{xcolor}
\usepackage{bm}
\usepackage{subfigure}

\providecommand{\etal}		{\emph{et al\@.}\xspace}
\providecommand{\ie}		{\emph{i.e\@.}\xspace}

\providecommand{\myurl}[1][]	{\texttt{web.eecs.umich.edu/$\sim$fessler#1}\xspace}

\providecommand{\onweb}[1]	{Available from \myurl.}

\long\def\comment#1{}

\providecommand{\bcent}		{\begin{center}}
\providecommand{\ecent}		{\end{center}}
\providecommand{\benum}		{\begin{enumerate}}
\providecommand{\eenum}		{\end{enumerate}}
\providecommand{\bitem}		{\begin{itemize}}
\providecommand{\eitem}		{\end{itemize}}
\providecommand{\bvers}		{\begin{verse}}
\providecommand{\evers}		{\end{verse}}
\providecommand{\btab}		{\begin{tabbing}}	
\providecommand{\etab}		{\end{tabbing}}

\newcounter{blist}
\providecommand{\blistmark}	{\makebox[0pt]{$\bullet$}}
\providecommand{\blistitemsep}	{0pt}
\providecommand{\blist}[1][]	{%
\begin{list}{\blistmark}{%
\usecounter{blist}%
\setlength{\itemsep}{\blistitemsep}%
\setlength{\parsep}{0pt}%
\setlength{\parskip}{0pt}%
\setlength{\partopsep}{0pt}%
\setlength{\topsep}{0pt}%
\setlength{\leftmargin}{1.2em}%
\setlength{\labelsep}{0.5\leftmargin}
\setlength{\labelwidth}{0em}%
#1}
}
\providecommand{\elist}		{\end{list}}

\providecommand{\blistitemsep}	{0pt}
\providecommand{\bjfenum}[1][]	{%
\begin{list}{\bcolor{\arabic{blist}.} }{%
\usecounter{blist}%
\setlength{\itemsep}{\blistitemsep}%
\setlength{\parsep}{0pt}%
\setlength{\parskip}{0pt}%
\setlength{\partopsep}{0pt}%
\setlength{\topsep}{0pt}%
\setlength{\leftmargin}{0.0em}%
\setlength{\labelsep}{1.0\leftmargin}
\setlength{\labelwidth}{0pt}%
#1}
}

\newcounter{blistAlph}
\providecommand{\blistAlph}[1][]
{\begin{list}{\makebox[0pt][l]{\Alph{blistAlph}.}}{%
\usecounter{blistAlph}%
\setlength{\itemsep}{0pt}\setlength{\parsep}{0pt}%
\setlength{\parskip}{0pt}\setlength{\partopsep}{0pt}%
\setlength{\topsep}{0pt}%
\setlength{\leftmargin}{1.2em}%
\setlength{\labelsep}{1.0\leftmargin}
\setlength{\labelwidth}{0.0\leftmargin}#1}%
}

\newcounter{blistRoman}
\providecommand{\blistRoman}[1][]
{\begin{list}{\Roman{blistRoman}.}{%
\usecounter{blistRoman}%
\setlength{\itemsep}{0.5em}\setlength{\parsep}{0pt}%
\setlength{\parskip}{0pt}\setlength{\partopsep}{0pt}%
\setlength{\topsep}{0pt}%
\setlength{\leftmargin}{4em}%
\setlength{\labelsep}{0.4\leftmargin}
\setlength{\labelwidth}{0.6\leftmargin}#1}%
}

\usepackage{bbm} 
\providecommand{\jfbbm}[1]	{\xmath{\mathbbm{#1}}} 
\providecommand{\qed}[1][0pt]	{\hfill\raisebox{#1}{\inmath{\Box}}} 

\providecommand{\reals}		{\jfbbm{R}}

\providecommand{\defequ}	{\stackrel{\bigtriangleup}{=}}

\providecommand{\norm}[1]	{\xmath{\left\| #1 \right\|}}

\providecommand{\Inprod}[2]	{\xmath{\left\langle #1,\ #2 \right\rangle}}

\let\equivsave\equiv
\def\equiv{\xmath{\equivsave}}

\providecommand{\ba}[1]		{\left[ \begin{array}{#1}}
\providecommand{\ea}		{\end{array} \right]}
\providecommand{\be}		{\begin{equation}}
\providecommand{\ee}[1]		{\label{#1}\end{equation}}
\providecommand{\bea}		{\begin{eqnarray}}
\providecommand{\eea}[1]	{\label{#1}\end{eqnarray}}
\providecommand{\beas}		{\begin{eqnarray*}}
\providecommand{\eeas}		{\end{eqnarray*}}
\providecommand{\beals}[1][1]	{\begin{alignat*}{#1}}	
\providecommand{\eeals}		{\end{alignat*}}

\providecommand{\berr}[2]{
\bgroup
\renewcommand{\theequation}{#1}
\be
#2
\ee{e,#1}
\egroup
\ignorespaces
}

\providecommand{\bearr}[2]{
\bgroup
\renewcommand{\theequation}{#1}
\bea
#2
\eea{e,#1}
\egroup
\ignorespaces
}

\providecommand{\inmath}	{\ensuremath}
\providecommand{\xmath}[1]	{\inmath{#1}\xspace}
\providecommand{\bmath}[1]	{\xmath{\bm{#1}}}	

\providecommand{\paren}[1]	{\xmath{\left(#1\right)}}

\providecommand{\braces}[1]	{\xmath{\left\{#1\right\}}}

\providecommand{\xfrac}[2]	{\xmath{\frac{#1}{#2}}}

\providecommand{\Frac}[2]	{\xmath{{#1}/{#2}}}
\providecommand{\pfrac}[3][]	{\xmath{\!\paren{#1\frac{#2}{#3}}}}

\providecommand{\sinc}		{\Ofun{\mathrm{sinc}}}

\usepackage{jf-accent}
\usepackage{jf-fun}
\usepackage{jf-names}

\newcommand{\st} {\xmath{\text{s.t.}\:}}

\renewcommand{\defequ} {\triangleq}

\renewcommand{\defequ} {\triangleq}

\newcommand{\DT} {DT}
\newcommand{\FO} {FO}

\newcommand{\lr} {\xmath{LR^2}}

\newcommand{\bmlam} {\bmath{\lambda}}
\newcommand{\bmtau} {\bmath{{\tau}}}

\newcommand{\bmh} {\bmath{h}}
\newcommand{\bmg} {\bmath{g}}
\newcommand{\bmG} {\bmath{G}}
\newcommand{\bmdel} {\bmath{\delta}}
\newcommand{\bmu} {\bmath{u}}

\newcommand{\bme} {\bmath{e}}
\newcommand{\bmx} {\bmath{x}}
\newcommand{\bmy} {\bmath{y}}
\newcommand{\bmz} {\bmath{z}}
\newcommand{\bmnu} {\bmath{\nu}}

\newcommand{\bmS} {\bmath{S}}

\newcommand{\bmA} {\bmath{A}}

\newcommand{\bmD} {\bmath{D}}

\newcommand{\hopt} {\xmath{h_{\mathrm{opt}}}}
\newcommand{\hoptN} {\xmath{\hopt(N)}}

\newcommand{\fwoneGM} {\xmath{f_{\mathrm{1,GM}}(\bmx;h,N)}}
\newcommand{\fwoneGMz} {\xmath{f_{\mathrm{1,GM}}}}
\newcommand{\fwoneOGM} {\xmath{f_{\mathrm{1,OGM}}(\bmx;N)}}
\newcommand{\fwoneOGMz} {\xmath{f_{\mathrm{1,OGM}}}}
\newcommand{\fwoneOGMp} {\xmath{f_{\mathrm{1,OGM}'}}}
\newcommand{\fwtwo} {\xmath{f_{\mathrm{2}}(\bmx)}}
\newcommand{\fwtwoz} {\xmath{f_{\mathrm{2}}}}

\newcommand{\bmt} {\bmath{t}}

\newcommand{\hlam} {\xmath{\hat{\lambda}}}
\newcommand{\htau} {\xmath{\hat{\tau}}}
\newcommand{\hgam} {\xmath{\hat{\gamma}}}

\newcommand{\bmhlam} {\bmath{\hat{\lambda}}}
\newcommand{\bmhtau} {\bmath{{\hat{\tau}}}}

\newcommand{\cC} {\xmath{\mathcal{C}_L^{1,1}(\Reals^d)}}
\newcommand{\cF} {\xmath{\mathcal{F}_L(\Reals^d)}}

\newcommand{\cL} {\xmath{\mathcal{L}}}

\newcommand{\Reals} {\xmath{\mathbb{R}}}

\newcommand{\hki} {\xmath{h_{i,k}}}
\newcommand{\hkip} {\xmath{h_{i+1,k}}}

\newcommand{\hkj} {\xmath{h_{j,k}}}

\newcommand{\Hh} {\bmath{\hat{h}}}

\newcommand{\Hhki} {\xmath{\hat{h}_{i,k}}}
\newcommand{\Hhkip} {\xmath{\hat{h}_{i+1,k}}}
\newcommand{\Hhkj} {\xmath{\hat{h}_{j,k}}}

\newcommand{\Hhiip} {\xmath{\hat{h}_{i+1,i}}}
\newcommand{\Hhimi} {\xmath{\hat{h}_{i,i-1}}}

\journalname{JOTA}
\numberwithin{equation}{section}

\usepackage{fullpage}
\usepackage{hyperref}
\hypersetup{
        linktoc=page,
        colorlinks=true,
        linkcolor=blue,
        citecolor=red,
        filecolor=magenta,
        urlcolor=cyan
}

\begin{document}

\title{
On the convergence analysis of the optimized gradient method
\thanks{This research was supported in part by NIH grant U01 EB018753.}
}


\author{Donghwan Kim         \and
        Jeffrey A. Fessler 
}


\institute{Donghwan Kim \and Jeffrey A. Fessler \at
                Dept. of Electrical Engineering and Computer Science,
                University of Michigan, Ann Arbor, MI 48109 USA \\
                \email{kimdongh@umich.edu, fessler@umich.edu}           
}

\date{Date of current version: \today} 

\maketitle

\begin{abstract}
%
%
This paper considers the problem
of unconstrained minimization
of smooth convex functions
having Lipschitz continuous gradients
with known Lipschitz constant.
We recently proposed an optimized gradient method (OGM)
\cite{kim::ofo}
for this problem
and showed that it has a worst-case convergence bound
for the cost function decrease
that is twice as small
as that of Nesterov's fast gradient method (FGM)
\cite{nesterov:83:amf},
yet has a similarly efficient practical implementation.
Drori~\cite{drori:16:tei}
showed recently that OGM has optimal complexity
over the general class of first-order methods.
This optimality makes it important
to study fully the convergence properties of OGM.
The previous
worst-case convergence bound for OGM
was derived for only the last iterate
of a secondary sequence.
This paper provides an analytic convergence bound
for the primary sequence generated by OGM.
We then discuss additional convergence properties of OGM,
including the interesting fact that OGM
has two types of worst-case functions:
a piecewise affine-quadratic function and a quadratic function.
These results help complete the theory
of optimal first-order methods
for smooth convex minimization.
%
%
\keywords{First-order algorithms 
\and Optimized gradient methods
\and Convergence bound 
\and Smooth convex minimization \and Worst-case performance analysis}
\end{abstract}

\section{Introduction}
\label{intro}

Consider the unconstrained smooth convex minimization problem
\begin{align}
\min_{\bmx\in\Reals^d} \;&\; f(\bmx)
\label{eq:prob} \tag{M}
\end{align}
with the following three conditions:
\begin{itemize}[leftmargin=40pt]
\item
$f\;:\;\Reals^d\rightarrow\Reals$ is
a convex function of the type \cC,
\ie, continuously differentiable with Lipschitz continuous gradient:
\begin{align*}
||\nabla f(\bmx) - \nabla f(\bmy)|| \le L||\bmx-\bmy||, \quad \forall \bmx, \bmy\in\Reals^d
,\end{align*}
where $L > 0$ is the Lipschitz constant.
\item
The optimal set $X_*(f)=\argmin{\bmx\in\Reals^d} f(\bmx)$ is nonempty,
\ie, problem~\eqref{eq:prob} is solvable.
\item
The distance between the initial point $\bmx_0$
and an optimal solution $\bmx_* \in X_*(f)$ is bounded by $R>0$,
\ie, $||\bmx_0 - \bmx_*|| \le R$.
\end{itemize}
We use \cF to denote the class of functions that satisfy the above conditions hereafter.

For large-scale optimization problems of type~\eqref{eq:prob} 
that arise in various fields such as communications, 
machine learning and signal processing,
general first-order algorithms
that query only the cost function values and gradients
are attractive because of their mild dependence
on the problem dimension~\cite{cevher:14:cof}.
For simplicity,
we initially 
focus on
the class of \emph{fixed-step} first-order (\FO)
algorithms having the following form:

\begin{center}
\fbox{
\begin{minipage}[t]{0.47\textwidth}
\vspace{-10pt}
\begin{flalign}
&\text{\bf Algorithm Class \FO} & \nonumber \\
&\text{Input: } f\in \cF,\; \bmx_0\in\Reals^d. & \nonumber \\
&\text{For } i = 0,\ldots,N-1 & \nonumber \\
&\quad \bmx_{i+1} = \bmx_i - \frac{1}{L}\sum_{k=0}^i \hkip \nabla f(\bmx_k) & \label{eq:fo}
\end{flalign}
\end{minipage}
} \vspace{5pt}
\end{center}

\noindent
\FO~updates use weighted sums
of current and previous gradients $\{\nabla f(\bmx_k)\}_{k=0}^i$
with (pre-determined) step sizes $\{\hkip\}_{k=0}^i$ and the Lipschitz constant $L$.
Class \FO~includes the (fixed-step) gradient method (GM), 
the heavy-ball method~\cite{polyak:64:smo},
Nesterov's fast gradient method (FGM)~\cite{nesterov:83:amf,nesterov:05:smo},
and the recently introduced optimized gradient method (OGM)~\cite{kim::ofo}.
Those four methods
have efficient recursive formulations rather than
directly using~\eqref{eq:fo}
that would require storing all previous gradients
and computing weighted summations every iteration.
Among class \FO, Nesterov's FGM
has been used widely,
since it achieves the optimal rate $O(\Frac{1}{N^2})$
for decreasing a cost function in $N$ iterations~\cite{nesterov:04},
and has two efficient forms as shown below for smooth convex problems.

\begin{center}
\fbox{
\begin{minipage}[t]{0.47\textwidth}
\vspace{-10pt}
\begin{flalign*}
&\text{\bf Algorithm FGM1~\cite{nesterov:83:amf}} & \\
&\text{Input: } f\in \cF,\; \bmx_0\in\Reals^d,\;
        \bmy_0 = \bmx_0,\; t_0 = 1. & \\
&\text{For } i = 0,\ldots,N-1 & \\
&\quad \bmy_{i+1} = \bmx_i - \frac{1}{L}\nabla f(\bmx_i) & \\
&\quad t_{i+1} = \frac{1+\sqrt{1+4t_i^2}}{2}, & \\
&\quad \bmx_{i+1} = \bmy_{i+1}
                + \frac{t_i - 1}{t_{i+1}}(\bmy_{i+1} - \bmy_i) 
\end{flalign*}
\end{minipage}\vline\hspace{8pt}%
\begin{minipage}[t]{0.47\textwidth}
\vspace{-10pt}
\begin{flalign*}
&\text{\bf Algorithm FGM2~\cite{nesterov:05:smo}} & \\
&\text{Input: } f\in \cF,\; \bmx_0\in\Reals^d,\;
        \bmy_0 = \bmx_0,\; t_0 = 1. & \\
&\text{For } i = 0,\ldots,N-1 & \\
&\quad \bmy_{i+1} = \bmx_i - \frac{1}{L}\nabla f(\bmx_i) & \\
&\quad \bmz_{i+1} = \bmx_0 - \frac{1}{L}\sum_{k=0}^it_k \nabla f(\bmx_k) & \\
&\quad t_{i+1} = \frac{1+\sqrt{1+4t_i^2}}{2}, & \\
&\quad \bmx_{i+1} = \paren{1 - \frac{1}{t_{i+1}}}\bmy_{i+1}
                + \frac{1}{t_{i+1}}\bmz_{i+1} &
\end{flalign*}
\end{minipage}
} \vspace{5pt}
\end{center}

\noindent
Both FGM1 and FGM2 produce identical sequences $\{\bmy_i\}$ and $\{\bmx_i\}$,
where the primary sequence $\{\bmy_i\}$
satisfies the following convergence bound~\cite{nesterov:83:amf,nesterov:05:smo}
for any $1\le i\le N$:
\begin{align}
f(\bmy_i) - f(\bmx_*) \le \frac{LR^2}{2t_{i-1}^2}
        \le \frac{2LR^2}{(i+1)^2}
\label{eq:fgmp,conv}
.\end{align}
In~\cite{kim::ofo}, 
we showed that the secondary sequence $\{\bmx_i\}$ of FGM
satisfies the following convergence bound that is similar to~\eqref{eq:fgmp,conv}
for any $1 \le i\le N$:
\begin{align}
f(\bmx_i) - f(\bmx_*) \le \frac{LR^2}{2t_i^2}
        \le \frac{2LR^2}{(i+2)^2}
\label{eq:fgms,conv}
.\end{align}
Taylor~\etal~\cite{taylor::ssc} demonstrated
that the upper bounds~\eqref{eq:fgmp,conv} and~\eqref{eq:fgms,conv}
are only asymptotically tight.

When the large-scale condition ``$d \ge 2N+1$'' holds,
Nesterov~\cite{nesterov:04} 
showed that for any first-order method
generating $\bmx_N$ after $N$ iterations
there exists a function $\phi$ in \cF 
that satisfies the following lower bound:
\begin{align}
\frac{3L||\bmx_0 - \bmx_*||^2}{32(N+1)^2} 
	\le \phi(\bmx_N) - \phi(\bmx_*)
\label{eq:optbound0}
.\end{align}
Although FGM achieves the optimal rate $O(\Frac{1}{N^2})$, 
one can still seek algorithms that improve upon the constant factor
in~\eqref{eq:fgmp,conv} and~\eqref{eq:fgms,conv},
in light of the gap between 
the bounds~\eqref{eq:fgmp,conv},~\eqref{eq:fgms,conv} of FGM
and the lower complexity bound~\eqref{eq:optbound0}.
Building upon Drori and Teboulle (hereafter ``\DT'')'s approach~\cite{drori:14:pof} 
of seeking \FO~methods that are faster than Nesterov's FGM
(reviewed in Section~\ref{sec:pep,ogm,bound}),
we recently proposed following two efficient formulations of OGM~\cite{kim::ofo}.

\begin{center}
\fbox{
\begin{minipage}[t]{0.47\textwidth}
\vspace{-10pt}
\begin{flalign*}
&\text{\bf Algorithm OGM1} & \\
&\text{Input: } f\in \cF,\; \bmx_0\in\Reals^d,\;
        \bmy_0 = \bmx_0,\; \theta_0 = 1. & \\
&\text{For } i = 0,\ldots,N-1 & \\
&\quad \bmy_{i+1} = \bmx_i - \frac{1}{L}\nabla f(\bmx_i) & \\
&\quad \theta_{i+1}
                = \begin{cases}
                        \frac{1+\sqrt{1+4\theta_i^2}}{2}, & i \le N-2 \\
                        \frac{1 + \sqrt{1+8\theta_i^2}}{2}, & i = N-1
                        \end{cases} & \\
&\quad \bmx_{i+1} = \bmy_{i+1}
                + \frac{\theta_i - 1}{\theta_{i+1}}(\bmy_{i+1} - \bmy_i) & \\
&\hspace{120pt}		
                + \frac{\theta_i}{\theta_{i+1}}(\bmy_{i+1} - \bmx_i) &
\end{flalign*}
\end{minipage}\vline\hspace{8pt}%
\begin{minipage}[t]{0.47\textwidth}
\vspace{-10pt}
\begin{flalign*}
&\text{\bf Algorithm OGM2} & \\
&\text{Input: } f\in \cF,\; \bmx_0\in\Reals^d,\;
        \bmy_0 = \bmx_0,\; \theta_0 = 1. & \\
&\text{For } i = 0,\ldots,N-1 & \\
&\quad \bmy_{i+1} = \bmx_i - \frac{1}{L}\nabla f(\bmx_i) & \\
&\quad \bmz_{i+1} = \bmx_0 - \frac{1}{L}\sum_{k=0}^i2\theta_k \nabla f(\bmx_k) & \\
&\quad \theta_{i+1}
                = \begin{cases}
                        \frac{1+\sqrt{1+4\theta_i^2}}{2}, & i \le N-2 \\
                        \frac{1 + \sqrt{1+8\theta_i^2}}{2}, & i = N-1
                        \end{cases} & \\
&\quad \bmx_{i+1} = \paren{1 - \frac{1}{\theta_{i+1}}}\bmy_{i+1}
		+ \frac{1}{\theta_{i+1}}\bmz_{i+1} &
\end{flalign*}
\end{minipage}
} \vspace{5pt}
\end{center}

\noindent
OGM1 and OGM2 
have computational efficiency comparable to FGM1 and FGM2,
and produce identical primary sequence $\{\bmy_i\}$ 
and secondary sequence $\{\bmx_i\}$.
The last iterate $\bmx_N$ of OGM satisfies
the following analytical worst-case bound~\cite[Theorem 2]{kim::ofo}:
\begin{align}
f(\bmx_N) - f(\bmx_*) \le \frac{LR^2}{2\theta_N^2}
        \le \frac{LR^2}{(N+1)(N+1+\sqrt{2})}
\label{eq:ogm,conv}
,\end{align}
which is twice as small as those for FGM in
\eqref{eq:fgmp,conv} and~\eqref{eq:fgms,conv}.
Recently for the condition ``$d \ge N+1$'', 
Drori~\cite{drori:16:tei}
showed that for any first-order method
there exists a function $\psi$ in \cF that cannot be minimized
faster than the following lower bound:
\begin{align}
\frac{L||\bmx_0 - \bmx_*||^2}{2\theta_N^2}
        \le \psi(\bmx_N) - \psi(\bmx_*)
\label{eq:optbound1}
,\end{align}
where $\bmx_N$ is the $N$th iterate of any first-order method.
This lower complexity bound~\eqref{eq:optbound1} improves
on~\eqref{eq:optbound0},
and exactly matches the bound~\eqref{eq:ogm,conv} of OGM,
showing that OGM achieves the \emph{optimal} worst-case bound of the cost function
for first-order methods when $d \ge N+1$.
What is remarkable about Drori's result
is that OGM was derived
by optimizing over the class FO
having fixed step sizes,
leading to~\eqref{eq:ogm,conv},
whereas Drori's lower bound in~\eqref{eq:optbound1}
is for the general class of first-order methods where the step sizes are arbitrary.
It is interesting that OGM
with its fixed step sizes
is optimal over the apparently much broader class.

Because OGM has such optimality,
it is desirable to understand its properties thoroughly.  
For example,
analytical bounds for the primary sequence $\{\bmy_i\}$ 
of OGM have not been studied previously,
although numerical bounds were discussed
by Taylor \etal~\cite{taylor::ssc}.
This paper provides analytical bounds for
the primary sequence of OGM,
augmenting the convergence analysis of $\bmx_N$ of OGM
given in~\cite{kim::ofo}.
We also relate OGM to another version of Nesterov's accelerated first-order method
in~\cite{nesterov:13:gmf}
that has a similar formulation as OGM2.

In~\cite[Theorem 3]{kim::ofo}, we specified a worst-case function
for which OGM achieves the first upper bound in~\eqref{eq:ogm,conv} exactly.
The corresponding worst-case function 
is the following piecewise affine-quadratic function:
\begin{align}
\fwoneOGM = \begin{cases}
		\frac{LR}{\theta_N^2}||\bmx|| - \frac{LR^2}{2\theta_N^4},
			& \text{if } ||\bmx|| \ge \frac{R}{\theta_N^2}, \\
		\frac{L}{2}||\bmx||^2, & \text{otherwise,}
	\end{cases}
\label{eq:fwoneOGM}
\end{align}
where OGM iterates remain in the affine region with the same gradient value
(without overshooting) for all $N$ iterations.
Section~\ref{sec:two,worst} shows that 
a simple quadratic function is also a worst-case function for OGM,
and describes
why it is interesting that the optimal OGM 
has these two types of worst-case functions.

Section~\ref{sec:bg} reviews
\DT's Performance Estimation Problem (PEP) framework in~\cite{drori:14:pof}
that enables systematic worst-case performance analysis of optimization methods.
Section~\ref{sec:ogm,second} provides new convergence analysis
for the primary sequence of OGM.
Section~\ref{sec:two,worst} discusses 
the two types of worst-case functions for OGM,
and Section~\ref{sec:conc} concludes.

\section{Prior work: Performance Estimation Problem (PEP)}
\label{sec:bg}

Exploring the convergence performance of optimization methods
including class \FO~has a long history. 
\DT~\cite{drori:14:pof}
were the first to cast the analysis of the worst-case performance of optimization methods 
into an optimization problem called PEP,
reviewed in this section.
We also review how we developed OGM~\cite{kim::ofo} 
that is built upon \DT's PEP.

\subsection{Review of PEP}
\label{sec:pep}

To analyze the worst-case convergence behavior
of a method in class \FO~having given step sizes $\bmh = \{\hki\}_{0\le k<i\le N}$,
\DT's PEP~\cite{drori:14:pof}
bounds the decrease of the cost function after $N$ iterations as
\begin{align}
\mathcal{B}_{\mathrm{P}}(\bmh,N,d,L,R) \defequ\;
&\max_{\substack{f\in\cF, \\ \bmx_0,\cdots,\bmx_N\in\Reals^d, \\
                \bmx_*\in X_*(f)}}
        f(\bmx_N) - f(\bmx_*)
\label{eq:PEP} \tag{P} \\
&\st \; ||\bmx_0 - \bmx_*|| \le R, \quad \bmx_{i+1} = \bmx_i - \frac{1}{L}\sum_{k=0}^i \hkip \nabla f(\bmx_k),
        \quad i=0,\ldots,N-1, \nonumber 
\end{align}
for given dimension $d$,
Lipschitz constant $L$
and the distance $R$ between an initial point $\bmx_0$
and an optimal point $\bmx_* \in X_*(f)$.

Since problem~\eqref{eq:PEP} is difficult to solve,
\DT~\cite{drori:14:pof} introduced a series of relaxations.
Then the upper bound of the worst-case performance
was found numerically in~\cite{drori:14:pof}
by solving a relaxed PEP problem.
For some cases,
analytical worst-case bounds were revealed in~\cite{drori:14:pof,kim::ofo},
where some of those analytical bounds
were even found to be exact despite the relaxations.
On the other hand, Taylor \etal~\cite{taylor::ssc}
studied the numerical tight worst-case bound of~\eqref{eq:PEP}
by avoiding \DT's one relaxation step that is not guaranteed to be tight
and showing the tightness of the rest of \DT's relaxations in~\cite{drori:14:pof}
(for the condition ``$d\ge N+2$'').

To summarize recent PEP studies,
\DT~extended the PEP approach
for nonsmooth convex problems~\cite{drori::aov},
Drori's thesis~\cite{drori:14:phd} includes
an extension of PEP
to projected gradient methods for constrained smooth convex problems,
and Taylor \etal~\cite{taylor:15:ewc}
studied PEP for various first-order algorithms
for solving composite convex problems.
Similarly but using different relaxations of~\eqref{eq:PEP},
Lessard \etal~\cite{lessard:16:aad}
applied the Integral quadratic constraints to~\eqref{eq:PEP},
leading to simpler computation
but slightly looser convergence upper bounds.

The next two sections review relaxations of \DT's PEP
and an approach for optimizing the choice of \bmh for \FO~using PEP
in~\cite{drori:14:pof}.

\subsection{Review of \DT's relaxation on PEP}
\label{sec:pep,relax}

This section reviews relaxations introduced by \DT~to make~\eqref{eq:PEP} 
into a simpler semidefinite programming (SDP) problem.
\DT~first relax the functional constraint $f\in\cF$ by
a well-known property of the class of \cF functions 
in~\cite[Theorem 2.1.5]{nesterov:04} and then further relax as follows:
%
\begin{align}
\mathcal{B}_{\mathrm{P1}}(\bmh,N,d,L,R) \defequ\;
&\max_{\substack{\bmG\in\Reals^{(N+1)d}, \\ \bmdel\in\Reals^{N+1}}}
        LR^2\delta_N
        \label{eq:pPEP} \tag{P1} \\
        &\st \; 
		\frac{1}{2}||\bmg_{i-1} - \bmg_i||^2 \le \delta_{i-1} - \delta_i
			- \Inprod{\bmg_i}{\sum_{k=0}^{i-1}\hki\bmg_k},
                	\quad i=1,\ldots,N, \nonumber \\
	&\qquad\!	\frac{1}{2}||\bmg_i||^2 \le -\delta_i
			- \Inprod{\bmg_i}{\sum_{j=1}^i\sum_{k=0}^{j-1}\hkj\bmg_k + \bmnu},
                	\quad i=0,\ldots,N, \nonumber
\end{align}
for any given unit vector 
$\bmnu\in\reals^d$, 
where we denote
$\bmg_i \defequ \frac{1}{L||\bmx_0 - \bmx_*||}\nabla f(\bmx_i)$
and
$\delta_i \defequ \frac{1}{L||\bmx_0 - \bmx_*||^2}(f(\bmx_i) - f(\bmx_*))$
for $i=0,\ldots,N,*$,
and define
$\bmG = [\bmg_0 \cdots \bmg_N]^\top \in \Reals^{(N+1)\times d}$
and $\bmdel = [\delta_0 \cdots \delta_N]^\top \in \Reals^{N+1}$.

Maximizing relaxed problem~\eqref{eq:pPEP} is still difficult,
so \DT~\cite{drori:14:pof} use a duality approach
on~\eqref{eq:pPEP}.
Replacing $\max_{\bmG,\bmdel} LR^2\delta_N$
by $\min_{\bmG,\bmdel} -\delta_N$ for convenience,
the Lagrangian 
of the corresponding
constrained minimization problem~\eqref{eq:pPEP} 
with dual variables $\bmlam = (\lambda_1,\cdots,\lambda_N)^\top\in\Reals_+^N$
and $\bmtau = (\tau_0,\cdots,\tau_N)^\top\in\Reals_+^{N+1}$
becomes
\begin{align}
\cL(\bmG,\bmdel,\bmlam,\bmtau;\bmh)
       \defequ 	- \delta_N + \sum_{i=1}^N \lambda_i(\delta_i - \delta_{i-1})
                	+ \sum_{i=0}^N \tau_i\delta_i 
		+ \Tr{\bmG^\top\bmS(\bmh,\bmlam,\bmtau)\bmG
                	+ \bmnu\bmtau^\top\bmG} 
\label{eq:lag}
,\end{align}
where
\begin{align}
\begin{cases}
\bmS(\bmh,\bmlam,\bmtau) \defequ \sum_{i=1}^N \lambda_i\bmA_{i-1,i}(\bmh)
                + \sum_{i=0}^N\tau_i\bmD_i(\bmh), & \\
\bmA_{i-1,i}(\bmh) \defequ \frac{1}{2}(\bmu_{i-1} - \bmu_i)(\bmu_{i-1} - \bmu_i)^\top
        + \frac{1}{2}\sum_{k=0}^{i-1}
        \hki (\bmu_i\bmu_k^\top + \bmu_k\bmu_i^\top), & \\
\bmD_i(\bmh) \defequ \frac{1}{2}\bmu_i\bmu_i^\top + \frac{1}{2}\sum_{j=1}^i\sum_{k=
0}^{j-1}
                \hkj(\bmu_i\bmu_k^\top + \bmu_k\bmu_i^\top), &
\end{cases}
\label{eq:ABCD}
\end{align}
and $\bmu_i = \bme_{i+1}\in\reals^{N+1}$
is the $(i+1)$th standard basis vector.

Using further derivations of a duality approach 
on~\eqref{eq:lag} in~\cite{drori:14:pof},
the dual problem of~\eqref{eq:pPEP} becomes the following SDP problem:
\begin{align}
\mathcal{B}_{\mathrm{P}}(\bmh,N,d,L,R) 
\le &\mathcal{B}_{\mathrm{D}}(\bmh,N,L,R) \defequ
\min_{\substack{(\bmlam,\bmtau)\in\Lambda, \\ \gamma\in\Reals}}
        \braces{ \frac{1}{2}LR^2\gamma\;:\;
        \paren{\begin{array}{cc}
        \bmS(\bmh,\bmlam,\bmtau) & \frac{1}{2}\bmtau \\
        \frac{1}{2}\bmtau^\top & \frac{1}{2}\gamma
        \end{array}} \succeq 0
	}
        \label{eq:dPEP} \tag{D}
,\end{align}
where
\begin{align*}
\Lambda = \braces{(\bmlam,\bmtau)\in\reals_+^N\times\reals_+^{N+1}\;:\;
                \begin{array}{l}
		\tau_0 = \lambda_1, \quad \lambda_N + \tau_N = 1, \\
                \lambda_i - \lambda_{i+1} + \tau_i = 0, \; i=1,\ldots,N-1, 
		\end{array}
        }
.\end{align*}
Then, for given \bmh, the bound $\mathcal{B}_{\mathrm{D}}(\bmh,N,L,R)$
(that is not guaranteed to be tight)
can be numerically computed using any SDP solver,
while analytical upper bounds  
$\mathcal{B}_{\mathrm{D}}(\bmh,N,L,R)$ 
for some choices of \bmh were found in
\cite{drori:14:pof,kim::ofo}.
Section~\ref{sec:ogm,second} finds a new analytical upper bound
for a modified version of $\mathcal{B}_{\mathrm{D}}$.

\subsection{Review of optimizing the step sizes using PEP}
\label{sec:pep,ogm,bound}

In addition to finding upper bounds for given \FO~methods,
\DT~\cite{drori:14:pof} searched for the best \FO~methods
with respect to the worst-case performance.
Ideally one would like to optimize \bmh
over problem~\eqref{eq:PEP}:
\begin{align}
\Hh_{\mathrm{P}} \defequ
\argmin{\bmh\in\Reals^{\Frac{N(N+1)}{2}}}
\mathcal{B}_{\mathrm{P}}(\bmh,N,d,L,R)
\label{eq:HP} \tag{HP}
.\end{align}
However, optimizing~\eqref{eq:HP} directly seems impractical,
so \DT~minimized the dual problem~\eqref{eq:dPEP} using a SDP solver
over the coefficients \bmh as
\begin{align}
\Hh_{\mathrm{D}} \defequ
\argmin{\bmh\in\Reals^{\Frac{N(N+1)}{2}}}
\mathcal{B}_{\mathrm{D}}(\bmh,N,L,R)
\label{eq:HD} \tag{HD}
.\end{align}
Due to relaxations, the computed $\Hh_{\mathrm{D}}$ 
is not guaranteed to be optimal for problem~\eqref{eq:HP}.
Nevertheless,
we show in~\cite{kim::ofo} that
solving~\eqref{eq:HD} leads to
an algorithm (OGM) having a convergence bound
that is twice as small as that of FGM.
Interestingly, OGM is optimal 
among first-order methods with $d\ge N+1$~\cite{drori:16:tei}, \ie,
$\Hh_{\mathrm{D}}$
is a solution of both~\eqref{eq:HP} and~\eqref{eq:HD}
for $d\ge N+1$.
An optimal point $(\Hh,\bmhlam,\bmhtau,\hgam)$ of~\eqref{eq:HD} 
is given in~\cite[Lemma 4 and Proposition 3]{kim::ofo} as follows:
\begin{align}
\Hhkip
&= \begin{cases}
        \frac{\theta_i - 1}{\theta_{i+1}} \Hhki,
        & k = 0,\ldots,i-2, \\
        \frac{\theta_i - 1}{\theta_{i+1}}(\Hhimi - 1),
        & k = i - 1, \\
        1 + \frac{2\theta_i - 1}{\theta_{i+1}},
        & k = i,
        \end{cases}
\label{eq:ogm1,h} \\
&= \begin{cases}
        \frac{1}{\theta_{i+1}}\paren{2\theta_k - \sum_{j=k+1}^i \Hhkj},
        & k = 0,\ldots,i-1, \\
        1 + \frac{2\theta_i - 1}{\theta_{i+1}},
        & k = i,
\end{cases} 
\label{eq:ogm2,h} \\
\hat{\lambda}_i &= 
        \frac{2\theta_{i-1}^2}{\theta_N^2}, \quad i=1,\ldots,N, \quad
\hat{\tau}_i = \begin{cases}
	\frac{2\theta_i}{\theta_N^2}, & i=0,\ldots,N-1, \\
	\frac{1}{\theta_N}, & i = N,
	\end{cases} \quad
\hgam = \frac{1}{\theta_N^2}.
\end{align}
Thus both OGM1 and OGM2 satisfy the convergence bound
\eqref{eq:ogm,conv}~\cite[Theorem 2, Propositions 4 and 5]{kim::ofo}.

\section{New convergence analysis for the primary sequence of OGM}
\label{sec:ogm,second}

\subsection{Relaxed PEP for the primary sequence of OGM}
\label{sec:pep,ogm,second}

This section applies PEP to an iterate $\bmy_N$
of the following class of fixed-step first-order methods (\FO$'$),
complementing the worst-case performance of $\bmx_N$ in the previous section.

\begin{center}
\fbox{
\begin{minipage}[t]{0.48\textwidth}
\vspace{-10pt}
\begin{flalign}
&\text{\bf Algorithm Class \FO$'$} & \nonumber \\
&\text{Input: } f\in \cF,\; \bmx_0\in\Reals^d,\; \bmy_0 = \bmx_0. & \nonumber \\
&\text{For } i = 0,\ldots,N & \nonumber \\
&\quad \bmy_{i+1} = \bmx_i - \frac{1}{L}\nabla f(\bmx_i) & \nonumber \\
&\quad \bmx_{i+1} = \bmx_i - \frac{1}{L}\sum_{k=0}^i \hkip \nabla f(\bmx_k). & \nonumber
\end{flalign}
\end{minipage}
} \vspace{5pt}
\end{center}

We first replace $f(\bmx_N) - f(\bmx_*)$ in~\eqref{eq:PEP} 
by $f(\bmy_{N+1}) - f(\bmx_*)$ as follows:
\begin{align}
\mathcal{B}_{\mathrm{P'}}(\bmh,N,d,L,R) \defequ\;
&\max_{\substack{f\in\cF, \\ \bmx_0,\cdots,\bmx_N,\bmy_{N+1}\in\Reals^d, \\
                \bmx_*\in X_*(f)}}
        f(\bmy_{N+1}) - f(\bmx_*)
\label{eq:PEP_} \tag{P$'$} \\
&\st \; 
	||\bmx_0 - \bmx_*|| \le R, \quad \bmy_{N+1} = \bmx_N - \frac{1}{L}\nabla f(\bmx_N), 
		\nonumber \\
&\qquad\!	\bmx_{i+1} = \bmx_i - \frac{1}{L}\sum_{k=0}^i \hkip \nabla f(\bmx_k),
		\quad i=0,\ldots,N-1. \nonumber
\end{align}
We could directly repeat relaxations on~\eqref{eq:PEP_} 
as reviewed in Section~\ref{sec:pep,relax},
but we found it difficult to solve a such relaxed problem of~\eqref{eq:PEP_} analytically.
Instead, we use the following inequality~\cite{nesterov:04}:
\begin{align}
f\paren{\bmx - \frac{1}{L}\nabla f(\bmx)} \le f(\bmx) - \frac{1}{2L}\norm{\nabla f(\bmx)}^2,
\quad \forall \bmx\in\Reals^d
\label{eq:ineqrelax}
.\end{align}
to relax~\eqref{eq:PEP_}, 
leading to the following bound:
\begin{align}
\mathcal{B}_{\mathrm{P1}'}(\bmh,N,d,L,R) \defequ\;
&\max_{\substack{f\in\cF, \\ \bmx_0,\cdots,\bmx_N\in\Reals^d, \\
                \bmx_*\in X_*(f)}}
        f(\bmx_N) - \frac{1}{2L}||\nabla f(\bmx_N)||^2 - f(\bmx_*)
\label{eq:PEP1_} \tag{P1$'$} \\
&\st \; ||\bmx_0 - \bmx_*|| \le R, \nonumber \\
     &\qquad\!   \bmx_{i+1} = \bmx_i - \frac{1}{L}\sum_{k=0}^i \hkip \nabla f(\bmx_k),
     \quad i=0,\ldots,N-1. \nonumber
\end{align}
This bound has an additional term $- \frac{1}{2L}||\nabla f(\bmx_N)||^2$
compared to~\eqref{eq:PEP}.
We later show that the increase of the worst-case upper bound 
due to this strict relaxation step using~\eqref{eq:ineqrelax}
is negligible asymptotically. 

Similar to relaxing from~\eqref{eq:PEP} to~\eqref{eq:pPEP} in Section~\ref{sec:pep,relax}, 
we relax~\eqref{eq:PEP1_} to the following bound:
\begin{align}
\mathcal{B}_{\mathrm{P2}'}(\bmh,N,d,L,R) \defequ\;
&\max_{\substack{\bmG\in\Reals^{(N+1)d}, \\ \bmdel\in\Reals^{N+1}}}
        LR^2\paren{\delta_N - \frac{1}{2}||\bmg_N||^2}
        \label{eq:pPEP_} \tag{P2$'$} \\
        &\st \;
		\frac{1}{2}||\bmg_{i-1} - \bmg_i||^2 \le \delta_{i-1} - \delta_i
                	- \Inprod{\bmg_i}{\sum_{k=0}^{i-1}\hki\bmg_k},
                	\quad i=1,\ldots,N, \nonumber \\
	&\qquad\! \frac{1}{2}||\bmg_i||^2 \le -\delta_i
                	- \Inprod{\bmg_i}{\sum_{j=1}^i\sum_{k=0}^{j-1}\hkj\bmg_k + \bmnu},
                	\quad i=0,\ldots,N, \nonumber
\end{align}
for any given unit vector $\bmnu\in\Reals^d$.
Then, as in Section~\ref{sec:pep,relax}
and~\cite{drori:14:pof,kim::ofo}, one can show that the dual problem of~\eqref{eq:pPEP_} 
is the following SDP problem
\begin{align}
\mathcal{B}_{\mathrm{P}'}(\bmh,N,d,L,R) 
\le &\mathcal{B}_{\mathrm{D}'}(\bmh,N,L,R) \defequ
\min_{\substack{(\bmlam,\bmtau)\in\Lambda, \\ \gamma\in\Reals}}
        \braces{ \frac{1}{2}LR^2\gamma\;:\;
        \paren{\begin{array}{cc}
        \bmS(\bmh,\bmlam,\bmtau) + \frac{1}{2}\bmu_N\bmu_N^\top & \frac{1}{2}\bmtau \\
        \frac{1}{2}\bmtau^\top & \frac{1}{2}\gamma
        \end{array}} \succeq 0}
	\label{eq:dPEP_} \tag{D$'$}
,\end{align}
by considering that the Lagrangian of~\eqref{eq:pPEP_} becomes
\begin{align}
\cL'(\bmG,\bmdel,\bmlam,\bmtau;\bmh) 
       \defequ  &- \delta_N + \sum_{i=1}^N \lambda_i(\delta_i - \delta_{i-1})
                        + \sum_{i=0}^N \tau_i\delta_i, 
                + \Tr{\bmG^\top\paren{\bmS(\bmh,\bmlam,\bmtau) 
			+ \frac{1}{2}\bmu_N\bmu_N^\top}\bmG
                        + \bmnu\bmtau^\top\bmG}
\label{eq:lag_}
\end{align}
when we replace $\max_{\bmG,\bmdel}LR^2\paren{\del_N - \frac{1}{2}||\bmg_N||^2}$
in~\eqref{eq:pPEP_}
by $\min_{\bmG,\bmdel}\paren{-\del_N + \frac{1}{2}||\bmg_N||^2}$ 
for simplicity as we did for~\eqref{eq:pPEP} and~\eqref{eq:lag}.
The formulation~\eqref{eq:lag_} is similar to~\eqref{eq:lag},
except the term $\frac{1}{2}\bmu_N\bmu_N^\top$.
The derivation of~\eqref{eq:dPEP_} and~\eqref{eq:lag_} is omitted here,
since it is almost identical to
the derivation of~\eqref{eq:dPEP} and~\eqref{eq:lag}
in~\cite{drori:14:pof,kim::ofo}.

\subsection{Convergence analysis for the primary sequence of OGM}
\label{sec:conv,ogm,second}

To find an upper bound for~\eqref{eq:dPEP_},
it suffices to specify a feasible point.

\begin{lemma}
\label{lem:hd}
The following choice of $(\Hh',\bmhlam',\bmhtau',\hgam')$
is a feasible point of~\eqref{eq:dPEP_}:
\begin{align}
\Hhkip'
&= \begin{cases}
        \frac{t_i - 1}{t_{i+1}} \Hhki',
        & k = 0,\ldots,i-2, \\
        \frac{t_i - 1}{t_{i+1}}(\Hhimi' - 1),
        & k = i - 1, \\
        1 + \frac{2t_i - 1}{t_{i+1}},
        & k = i,
        \end{cases}
\label{eq:ogm1_,h} \\
&= \begin{cases}
        \frac{1}{t_{i+1}}\paren{2t_k - \sum_{j=k+1}^i \Hhkj'},
        & k = 0,\ldots,i-1, \\
        1 + \frac{2t_i - 1}{t_{i+1}},
        & k = i,
\end{cases}
\label{eq:ogm2_,h} \\
\hlam_i' &=
        \frac{t_{i-1}^2}{t_N^2}, \quad i=1,\ldots,N, \quad
\htau_i' =
	\frac{t_i}{t_N^2}, \quad i=0,\ldots,N, \quad
\hgam' =
        \frac{1}{2t_N^2}.
\label{eq:hdvar_}
\end{align}
\begin{proof}
The equivalency between~\eqref{eq:ogm1_,h} and~\eqref{eq:ogm2_,h}
follows from~\cite[Proposition 3]{kim::ofo}.
Also, it is obvious that $(\bmhlam',\bmhtau')\in\Lambda$ using $t_i^2 = \sum_{k=0}^i t_k$.

We next rewrite $\bmS(\Hh',\bmhlam',\bmhtau')$
to show that
the choice $(\Hh',\bmhlam',\bmhtau',\hgam')$
satisfies the positive semidefinite condition in~\eqref{eq:dPEP_}.
For any \bmh and $(\bmlam,\bmtau)\in\Lambda$,
the $(i,k)$th entry of
the symmetric matrix $\bmS(\bmh,\bmlam,\bmtau)$ in~\eqref{eq:ABCD}
can be written as
\begin{align}
S_{i,k}(\bmh,\bmlam,\bmtau)
&= \begin{cases}
        \frac{1}{2}\Big((\lambda_i + \tau_i)\hki
                + \tau_i\sum_{j=k+1}^{i-1} \hkj\Big),
        & i=2,\ldots,N,\;k =0,\ldots, i - 2, \\
        \frac{1}{2}\paren{(\lambda_i + \tau_i)\hki - \lambda_i},
        & i=1,\ldots,N,\;k = i-1, \\
        \lambda_{i+1},
        & i=0,\ldots,N-1,\;k = i, \\
        \frac{1}{2},
        & i=N,\;k = i.
        \end{cases}
\label{eq:S2}
\end{align}
Inserting $\Hh'$, $\bmhlam'$ and $\bmhtau'$
into~\eqref{eq:S2}, we get
\begin{align*}
S_{i,k}(\Hh',\bmhlam',\bmhtau') + \frac{1}{2}\bmu_N\bmu_N^\top
= &\begin{cases}
        \frac{1}{2}\Big(\frac{t_{i}^2}{t_N^2}
                        \frac{1}{t_{i}}\paren{2t_k - \sum_{j=k+1}^{i-1}\Hhkj'}
                        + \frac{t_{i}}{t_N^2}\sum_{j=k+1}^{i-1}\Hhkj'\Big),
        & i=2,\ldots,N,\;k =0,\ldots,i - 2, \\
        \frac{1}{2}\paren{\frac{t_{i}^2}{t_N^2}
                \paren{1 + \frac{2t_{i-1} - 1}{t_{i}}} - \frac{t_{i-1}^2}{t_N^2}},
        & i=1,\ldots,N,\;k = i-1, \\
        \frac{t_i^2}{t_N^2},
        & i=0,\ldots,N-1,\;k = i, \\
        1,
        & i=N,\;k = i,
        \end{cases} \\
= &\frac{t_it_k}{t_N^2}
\end{align*}
where
the second equality uses $t_i^2 - t_i - t_{i-1}^2 = 0$.

Finally, using $\hgam'$, we have
\begin{align*}
\paren{\begin{array}{cc}
        \bmS(\Hh',\bmhlam',\bmhtau') + \frac{1}{2}\bmu_N\bmu_N^\top & \frac{1}{2}\bmhtau' \\
        \frac{1}{2}\bmhtau'^\top & \frac{1}{2}\hgam'
\end{array}}
= \paren{\begin{array}{cc} \frac{1}{t_N^2}\bmt\,\bmt^\top
                & \frac{1}{2t_N^2}\bmt \\
                \frac{1}{2t_N^2}\bmt^\top
                & \frac{1}{4t_N^2} \end{array}}
= \frac{1}{t_N^2}\paren{\begin{array}{c} \bmt \\ \frac{1}{2} \end{array}}
                \paren{\begin{array}{c} \bmt \\ \frac{1}{2} \end{array}}^\top
\succeq 0
,\end{align*}
where $\bmt = (t_0, \cdots, t_N)^\top$.
\qed
\end{proof}
\end{lemma}

Since \Hh~\eqref{eq:ogm1,h} 
and $\Hh'$~\eqref{eq:ogm1_,h} 
are identical except for the last iteration,
the intermediate iterates $\{\bmx_i\}_{i=0}^{N-1}$
of \FO~with both \Hh and $\Hh'$ are equivalent.
We can also easily notice that the sequence $\{\bmy_i\}_{i=0}^N$
of \FO$'$ with both \Hh and $\Hh'$ are also identical,
implying that both the primary sequence $\{\bmy_i\}$ of OGM
and \FO$'$ with $\Hh'$ are equivalent. 

Using Lemma~\ref{lem:hd},
the following theorem provides
an analytical convergence bound for
the primary sequence $\{\bmy_i\}$ of OGM.

\begin{theorem}
\label{thm:ogm_,conv}
Let $f\in\cF$
and let $\bmy_0,\cdots,\bmy_N \in \Reals^d$ be generated by OGM1 and OGM2. 
Then for $1\le i\le N$,
the primary sequence for OGM satisfies:
\begin{align}
f(\bmy_i) - f(\bmx_*) \le \frac{LR^2}{4t_{i-1}^2}
\le \frac{LR^2}{(i+1)^2}
\label{eq:ogm_,conv}
.\end{align}
\begin{proof}
The sequence $\{\bmy_i\}_{i=0}^N$ generated by \FO$'$ with $\Hh'$
is equivalent to that of OGM1 and OGM2~\cite[Propositions 4 and 5]{kim::ofo}.

Using $\hgam'$~\eqref{eq:hdvar_}
and $t_i^2 \ge \frac{(i+2)^2}{4}$,
we have
\begin{align}
f(\bmy_N) - f(\bmx_*) \le
\mathcal{B}_{\mathrm{D'}}(\Hh',N-1,L,R)
= \frac{1}{2}LR^2\hgam'
= \frac{LR^2}{4t_{N-1}^2}
\le \frac{LR^2}{(N+1)^2}
\label{eq:ogm_,conv0}
,\end{align}
based on Lemma~\ref{lem:hd}.
Since the primary sequence $\{\bmy_i\}_{i=0}^N$ of OGM1 and OGM2 
does not depend on a given $N$,
we can extend~\eqref{eq:ogm_,conv0}
for all $1\le i\le N$.
\qed
\end{proof}
\end{theorem}

Due to a strict relaxation leading to~\eqref{eq:PEP1_},
we cannot guarantee that the bound~\eqref{eq:ogm_,conv} is tight.
However, 
the next proposition shows that bound~\eqref{eq:ogm_,conv} is asymptotically tight
by specifying one particular worst-case function
that was conjectured by Taylor~\etal~\cite[Conjecture 4]{taylor::ssc}.

\begin{proposition}
\label{prop:ogmp,worst}
For the following function in \cF:
\begin{align}
\fwoneOGMp(\bmx;N) = \begin{cases}
        \frac{LR}{2t_{N-1}^2+1}||\bmx|| - \frac{LR^2}{2(2t_{N-1}^2+1)^2},
                & \text{if } ||\bmx||\ge\frac{R}{2t_{N-1}^2+1}, \\
        \frac{L}{2}||\bmx||^2,
                & \text{otherwise},
	\end{cases}
\label{eq:fwoneOGMp}
\end{align}
the iterate $\bmy_N$ generated by OGM1 and OGM2
provides the following lower bound:
\begin{align}
\frac{LR^2}{4t_{N-1}^2+2} 
	= \fwoneOGMp(\bmy_N;N) - \fwoneOGMp(\bmx_*;N)
	\le \max_{\substack{f\in\cF, \\ \bmx_* \in X_*(f)}} f(\bmy_N) - f(\bmx_*)
\label{eq:ogm_,low}
.\end{align}
\begin{proof}
Starting from $\bmx_0 = R\bmnu$, where \bmnu is a unit vector,
and using the following property of the coefficients $\Hh'$
\cite[Equation (8.2)]{kim::ofo}:
\begin{align}
\sum_{j=1}^i\sum_{k=0}^{j-1}\Hhkj' = t_i^2 - 1, \quad i=1,\ldots,N,
\end{align}
the primary iterates of OGM1 and OGM2 are as follows
\begin{align*}
\bmy_i &= \bmx_{i-1} - \frac{1}{L}\nabla\fwoneOGMp(\bmx_{i-1};N) 
	= \bmx_0 - \frac{1}{L}\sum_{j=1}^{i-1}\sum_{k=0}^{j-1}\Hhkj'\nabla\fwoneOGMp(\bmx_k;N)
		- \frac{1}{L}\nabla\fwoneOGMp(\bmx_{i-1};N) \\
	&= \paren{1 - \frac{t_{i-1}^2}{2t_{N-1}^2+1}}R\bmnu, \quad i=1,\ldots,N,
\end{align*}
where the corresponding sequence $\{\bmx_0,\cdots,\bmx_{N-1},\bmy_1,\cdots,\bmy_N\}$
stays in the affine region of the function $\fwoneOGMp(\bmx;N)$
with the same gradient value:
\begin{align*}
\nabla\fwoneOGMp(\bmx_i;N) = \nabla\fwoneOGMp(\bmy_{i+1};N) 
	= \frac{LR}{2t_{N-1}^2+1}\bmnu, \quad i=0,\ldots,N-1.
\end{align*}
Therefore, after $N$ iterations of OGM1 and OGM2, we have
\begin{align*}
\fwoneOGMp(\bmy_N;N) - \fwoneOGMp(\bmx_*;N) 
	= \fwoneOGMp\paren{\frac{t_{N-1}^2+1}{2t_{N-1}^2+1}R\bmnu;N}
        = \frac{LR^2}{2(2t_{N-1}^2+1)}
,\end{align*}
exactly matching the lower bound~\eqref{eq:ogm_,low}.
\qed
\end{proof}
\end{proposition}

The lower bound~\eqref{eq:ogm_,low} matches
the tight numerical worst-case bound in~\cite{taylor::ssc}
(see Table~\ref{tab:bound}).
While Taylor~\etal~\cite{taylor::ssc} provide
numerical evidence about the tight bound of the primary sequence of OGM,
our~\eqref{eq:ogm_,low}
provides an analytical bound that 
suffices for asymptotically tight worst-case analysis.

\subsection{New formulations of OGM}

Using~\cite[Propositions 4 and 5]{kim::ofo},
Algorithm \FO$'$ with the coefficients $\Hh'$~\eqref{eq:ogm1_,h} and~\eqref{eq:ogm2_,h}
can be implemented efficiently as follows:

\begin{center}
\fbox{
\begin{minipage}[t]{0.47\textwidth}
\vspace{-10pt}
\begin{flalign*}
&\text{\bf Algorithm OGM1$'$} & \\
&\text{Input: } f\in \cF,\; \bmx_0\in\Reals^d,\;
        \bmy_0 = \bmx_0,\; t_0 = 1. & \\
&\text{For } i = 0,\ldots,N-1 & \\
&\quad \bmy_{i+1} = \bmx_i - \frac{1}{L}\nabla f(\bmx_i) & \\
&\quad t_{i+1} = \frac{1+\sqrt{1+4t_i^2}}{2} & \\
&\quad \bmx_{i+1} = \bmy_{i+1}
                + \frac{t_i - 1}{t_{i+1}}(\bmy_{i+1} - \bmy_i) & \\
&\hspace{120pt}
		+ \frac{t_i}{t_{i+1}}(\bmy_{i+1} - \bmx_i) &
\end{flalign*}
\end{minipage}\vline\hspace{8pt}%
\begin{minipage}[t]{0.47\textwidth}
\vspace{-10pt}
\begin{flalign*}
&\text{\bf Algorithm OGM2$'$} & \\
&\text{Input: } f\in \cF,\; \bmx_0\in\Reals^d,\;
        \bmy_0 = \bmx_0,\;t_0=1. & \\
&\text{For } i = 0,\ldots,N-1 & \\
&\quad \bmy_{i+1} = \bmx_i - \frac{1}{L}\nabla f(\bmx_i) & \\
&\quad \bmz_{i+1} = \bmx_0 - \frac{1}{L}\sum_{k=0}^i2t_k\nabla f(\bmx_k) & \\
&\quad t_{i+1} = \frac{1+\sqrt{1+4t_i^2}}{2} & \\
&\quad \bmx_{i+1} = \paren{1 - \frac{1}{t_{i+1}}}\bmy_{i+1}
                                + \frac{1}{t_{i+1}}\bmz_{i+1} &
\end{flalign*}
\end{minipage}
} \vspace{5pt}
\end{center}
The OGM$'$ is very similar to OGM,
because it generates same primary and secondary sequence;
only the last iterate of the secondary sequence differs.
Therefore, the bound~\eqref{eq:ogm_,conv} applies to
the primary sequence $\{\bmy_i\}$ of both OGM and OGM$'$,
as summarized in the following corollary.

\begin{corollary}
\label{cor:ogm_,conv}
Let $f\in\cF$
and let $\bmy_0,\cdots, \bmy_N \in \Reals^d$ be generated by OGM1$'$ and OGM2$'$.
Then for $1\le i\le N$,
\begin{align}
f(\bmy_i) - f(\bmx_*) \le \frac{LR^2}{4t_{i-1}^2}
\le \frac{LR^2}{(i+1)^2}
\label{eq:ogm__,conv}
.\end{align}
\end{corollary}

\subsection{Comparing tight worst-case bounds of FGM, OGM and OGM$'$}
\label{sec:comp,worst}

While some analytical upper bounds of FGM, OGM and OGM$'$
such as \eqref{eq:fgmp,conv},~\eqref{eq:fgms,conv}
\eqref{eq:ogm,conv},~\eqref{eq:ogm_,conv}
and~\eqref{eq:ogm__,conv}
are available for comparison,
some of those are tight only asymptotically 
or some bounds for such algorithms are even unknown analytically.
Therefore, we used the code of Taylor~\etal~\cite{taylor::ssc} 
for tight (numerical) comparison of algorithms of interest
for some given $N$.
Table~\ref{tab:bound} provides tight numerical bounds of
the primary and secondary sequence of FGM, OGM and OGM$'$.
Interestingly, the worst-case performance
of the secondary sequence of OGM$'$ is
worse than that of FGM sequences,
whereas the primary sequence of OGM (and OGM$'$) is roughly twice better.

\begin{table}[h!]
\normalsize
\caption{Exact numerical bound of the last primary iterate $\bmy_N$
and the last secondary iterate $\bmx_N$ of FGM, OGM and OGM$'$}
\label{tab:bound}       
\begin{center}
\begin{tabular}{rlllll}
\hline\noalign{\smallskip}
$N$&FGM(primary)        &FGM(secondary)      &OGM(primary)        &OGM(secondary)      &OGM$'$(secondary) \\
\noalign{\smallskip}\hline\noalign{\smallskip}
 1 &\Frac{\lr}{6.00}    &\Frac{\lr}{6.00}    &\Frac{\lr}{6.00}    &\Frac{\lr}{8.00}    &\Frac{\lr}{5.24}    \\
 2 &\Frac{\lr}{10.00}   &\Frac{\lr}{11.13}   &\Frac{\lr}{12.47}   &\Frac{\lr}{16.16}   &\Frac{\lr}{9.62}    \\
 3 &\Frac{\lr}{15.13}   &\Frac{\lr}{17.35}   &\Frac{\lr}{21.25}   &\Frac{\lr}{26.53}   &\Frac{\lr}{15.12}   \\
 4 &\Frac{\lr}{21.35}   &\Frac{\lr}{24.66}   &\Frac{\lr}{32.25}   &\Frac{\lr}{39.09}   &\Frac{\lr}{21.71}   \\
 5 &\Frac{\lr}{28.66}   &\Frac{\lr}{33.03}   &\Frac{\lr}{45.42}   &\Frac{\lr}{53.80}   &\Frac{\lr}{29.38}   \\
10 &\Frac{\lr}{81.07}   &\Frac{\lr}{90.69}   &\Frac{\lr}{143.23}  &\Frac{\lr}{159.07}  &\Frac{\lr}{83.54}   \\
20 &\Frac{\lr}{263.65}  &\Frac{\lr}{283.55}  &\Frac{\lr}{494.68}  &\Frac{\lr}{525.09}  &\Frac{\lr}{269.56}  \\
40 &\Frac{\lr}{934.89}  &\Frac{\lr}{975.10}  &\Frac{\lr}{1810.08} &\Frac{\lr}{1869.22} &\Frac{\lr}{947.55}  \\
80 &\Frac{\lr}{3490.22} &\Frac{\lr}{3570.75} &\Frac{\lr}{6866.93} &\Frac{\lr}{6983.13} &\Frac{\lr}{3516.00} \\
\noalign{\smallskip}\hline
\end{tabular}
\end{center}
\end{table}

The following proposition uses a quadratic function
to define a lower bound on the worst-case performance 
of OGM1$'$ and OGM2$'$.

\begin{proposition}
For the following quadratic function in \cF:
\begin{align}
\fwtwo = \frac{L}{2}||\bmx||^2
\label{eq:fwtwo}
\end{align}
both OGM1$'$ and OGM2$'$ provide the following lower bound:
\begin{align}
\frac{LR^2}{2t_i^2}
	= \fwtwoz(\bmx_i) - \fwtwoz(\bmx_*)
	\le \max_{\substack{f\in\cF, \\ \bmx_* \in X_*(f)}} f(\bmx_i) - f(\bmx_*),
\label{eq:ogm,second,conv}
\end{align}
\begin{proof}
We use induction to show that
the following iterates:
\begin{align}
\bmx_i = (-1)^i\frac{1}{t_i}R\bmnu, \quad i=0,\ldots,N,
\label{eq:qpoint}
\end{align}
correspond to the iterates of OGM1$'$ and OGM2$'$ applied to \fwtwo.
Starting from $\bmx_0 = R\bmnu$, where \bmnu is a unit vector,
and assuming that~\eqref{eq:qpoint} holds for $i < N$,
we have
\begin{align*}
\bmx_{i+1} &= \bmx_i - \frac{1}{L}\sum_{k=0}^i\Hhkip'\nabla\fwtwoz(\bmx_k) \\
	&= \paren{\bmx_i - \frac{1}{L}\Hhiip'\nabla\fwtwoz(\bmx_i)}
                - \frac{1}{L}\sum_{k=0}^{i-1}\frac{t_i - 1}{t_{i+1}}
			\hat{h}_{i,k}'\nabla\fwtwoz(\bmx_k)
                + \frac{1}{L}\frac{t_i-1}{t_{i+1}}\nabla\fwtwoz(\bmx_{i-1}) \\
        &= \frac{1 - 2t_i}{t_{i+1}}\bmx_i
                + \frac{t_i-1}{t_{i+1}}(\bmx_i - \bmx_{i-1})
                + \frac{t_i-1}{t_{i+1}}\bmx_{i-1} 
	= - \frac{t_i}{t_{i+1}}\bmx_i \\
        &= (-1)^{i+1}\frac{1}{t_{i+1}}R\bmnu,
\end{align*}
where the second and third equalities use~\eqref{eq:fo} and~\eqref{eq:ogm1_,h}.
Therefore, we have
\begin{align*}
\fwtwoz(\bmx_N) - \fwtwoz(\bmx_*) 
= \fwtwoz\paren{(-1)^N\frac{1}{t_N}R\bmnu} = \frac{LR^2}{2t_N^2}
,\end{align*}
after $N$ iterations of OGM1$'$ and OGM2$'$, 
which is equivalent to the lower bound~\eqref{eq:ogm,second,conv}.
\qed
\end{proof}
\end{proposition}

Since the analytical lower bound~\eqref{eq:ogm,second,conv}
matches the numerical tight bound in Table~\ref{tab:bound},
we conjecture that the quadratic function \fwtwo
is the worst-case function for the secondary sequence of OGM$'$
and thus~\eqref{eq:ogm,second,conv} is the tight worst-case bound.
Whereas FGM has similar worst-case bounds
(and behavior as conjectured by Taylor~\etal~\cite[Conjectures 4 and 5]{taylor::ssc})
for both its primary and secondary sequence,
the two sequences of OGM$'$ (or intermediate iterates of OGM) 
have two different worst-case behaviors, 
as discussed further in Section~\ref{sec:two,worst,ogm}.

\subsection{Related work: 
Nesterov's accelerated first-order method in~\cite{nesterov:13:gmf}}

Interestingly, an algorithm in~\cite[Section 4]{nesterov:13:gmf}
is similar to OGM2$'$
and satisfies same convergence bound~\eqref{eq:ogm_,conv}
for the primary sequence $\{\bmy_i\}$,
which we call Nes13 in this paper for convenience.\footnote{
Nes13 was developed originally 
to deal with nonsmooth composite convex functions
with a line-search scheme~\cite[Section 4]{nesterov:13:gmf},
whereas the algorithm shown here 
is a simplified version of~\cite[Section 4]{nesterov:13:gmf}
for unconstrained smooth convex minimization~\eqref{eq:prob}
without a line-search.}

\begin{center}
\fbox{
\begin{minipage}[t]{0.47\textwidth}
\vspace{-10pt}
\begin{flalign*}
&\text{\bf Algorithm Nes13~\cite{nesterov:13:gmf}} & \\
&\text{Input: } f\in \cF,\; \bmx_0\in\Reals^d,\;
        \bmy_0 = \bmx_0,\;t_0=1. & \\
&\text{For } i = 0,1,\ldots & \\
&\quad \bmy_{i+1} = \bmx_i - \frac{1}{L}\nabla f(\bmx_i) & \\
&\quad \bmz_{i+1} = \bmx_0 - \frac{1}{L}\sum_{k=0}^i2t_k \nabla f(\bmy_{k+1}) & \\
&\quad t_{i+1} = \frac{1+\sqrt{1+4t_i^2}}{2} & \\
&\quad \bmx_{i+1} = \paren{1 - \frac{1}{t_{i+1}}}\bmy_{i+1}
                                + \frac{1}{t_{i+1}}\bmz_{i+1} &
\end{flalign*}
\end{minipage}
} \vspace{5pt}
\end{center}

\noindent
The only difference between OGM2$'$ and Nes13
is the gradient used for the update of $\bmz_i$.
While both algorithms achieve same bound~\eqref{eq:ogm_,conv},
Nes13 is less attractive in practice since
it requires computing gradients at two different points $\bmx_i$ and $\bmy_{i+1}$
at each $i$th iteration.


Similar to Proposition~\ref{prop:ogmp,worst},
the following proposition shows that the bound~\eqref{eq:ogm_,conv} is asymptotically tight
for Nes13.

\begin{proposition}
For the function $\fwoneOGMp(\bmx;N)$~\eqref{eq:fwoneOGMp} in \cF,
the iterate $\bmy_N$ generated by Nes13
achieves the lower bound~\eqref{eq:ogm_,low}.
\begin{proof}
See the proof of Proposition~\ref{prop:ogmp,worst}.
\qed
\end{proof}
\end{proposition}

\section{Two worst-case functions for an optimal fixed-step GM and OGM}
\label{sec:two,worst}

This section discusses two algorithms, an optimal fixed-step GM and OGM, 
in class \FO~that have a piecewise affine-quadratic function
and a quadratic function as two worst-case functions.
Considering that OGM is optimal among first-order methods (for $d\ge N+1$),
it is interesting
that OGM has two different types of worst-case functions,
because this property
resembles the (numerical) analysis 
of the optimal fixed-step GM in~\cite{taylor::ssc}
(reviewed below).

\subsection{Two worst-case functions for an optimal fixed-step GM}
\label{sec:two,worst,gm}

The following is GM with a constant step size $h$.

\begin{center}
\fbox{
\begin{minipage}[t]{0.47\textwidth}
\vspace{-10pt}
\begin{flalign}
&\text{\bf Algorithm GM} & \nonumber \\
&\text{Input: } f\in\cF,\; \bmx_0\in\Reals^d. & \nonumber \\
&\text{For } i = 0,\ldots,N & \nonumber \\
&\quad \bmx_{i+1} = \bmx_i - \frac{h}{L}\nabla f(\bmx_i) & \nonumber
\end{flalign}
\end{minipage}
} \vspace{5pt}
\end{center}

\noindent
For GM with $0<h<2$,
both \cite{drori:14:pof} and \cite{taylor::ssc} conjecture
the following tight convergence bound:
\begin{align}
f(\bmx_N) - f(\bmx_*) \le \frac{LR^2}{2}\max\paren{\frac{1}{2Nh+1}, (1-h)^{2N}}
\label{eq:gm,up}
.\end{align}
The proof of the bound~\eqref{eq:gm,up}
for $0< h\le 1$ is given in~\cite{drori:14:pof},
while the proof for $1<h<2$ is still unknown
but strong numerical evidence is given in~\cite{taylor::ssc}.
In other words,
at least one of the two functions specified below
is conjectured to be a worst-case
for GM with a constant step size $0<h<2$.
Such functions are
a piecewise affine-quadratic function
\begin{align}
\fwoneGM &= \begin{cases}
        \frac{LR}{2Nh+1}||\bmx|| - \frac{LR^2}{2(2Nh+1)^2},
                & \text{if } ||\bmx|| \ge \frac{R}{2Nh+1}, \\
        \frac{L}{2}||\bmx||^2, & \text{otherwise,}
\end{cases} \label{eq:fwoneh}
\end{align}
and a quadratic function
\fwtwo~\eqref{eq:fwtwo},
where \fwoneGM and \fwtwo contribute
to the factors $\frac{1}{2Nh+1}$ and $(1-h)^{2N}$ 
respectively in~\eqref{eq:gm,up}.
Here, \fwoneGM is a worst-case function
where the GM iterates approach the optimum slowly,
whereas \fwtwo is a worst-case function 
where the iterates overshoot the optimum.
(See Fig.~\ref{fig:gm,constant}.)

Assuming that the above conjecture for a fixed-step GM holds,
Taylor~\etal~\cite{taylor::ssc} searched (numerically) for
the optimal fixed-step size $0<h_{\mathrm{opt}}(N)<2$ for given $N$
that minimizes the bound~\eqref{eq:gm,up}: 
\begin{align}
h_{\mathrm{opt}}(N) \defequ \argmin{0<h<2} \max\paren{\frac{1}{2Nh+1}, (1-h)^{2N}}
\label{eq:hopt}
.\end{align}
GM with the step \hoptN
has two worst-case functions \fwoneGM and \fwtwo,
and must compromise between two extreme cases.
On the other hand,
the case $0<h<\hoptN$ has only \fwoneGM as the worst-case
and the case $\hoptN< h<2$ has only \fwtwo as the worst-case.
We believe this compromise 
is inherent to optimizing the worst-case performance of \FO~methods.
The next section shows that the optimal OGM 
also has this desirable property.

For the special case of $N=1$,
the optimal OGM 
reduces to GM with a fixed-step $h=1.5$,
and this confirms the conjecture in~\cite{taylor::ssc}
that the step $\hopt(1)=1.5$~\eqref{eq:hopt} is optimal
for a fixed-step GM with $N=1$.
However, proving the optimality of $\hopt(N)$~\eqref{eq:hopt}
for the fixed-step GM for $N>1$
is left as future work.

Fig.~\ref{fig:gm,constant} visualizes
the worst-case performance of GM with the optimal fixed-step \hoptN
for $N=2$ and $N=5$.
As discussed,
for the two worst-case function in Fig.~\ref{fig:gm,constant},
the final iterates
reach the same cost function value,
where the iterates approach the optimum slowly for \fwoneGM,
and overshoot for \fwtwo.

\begin{figure}[h!]
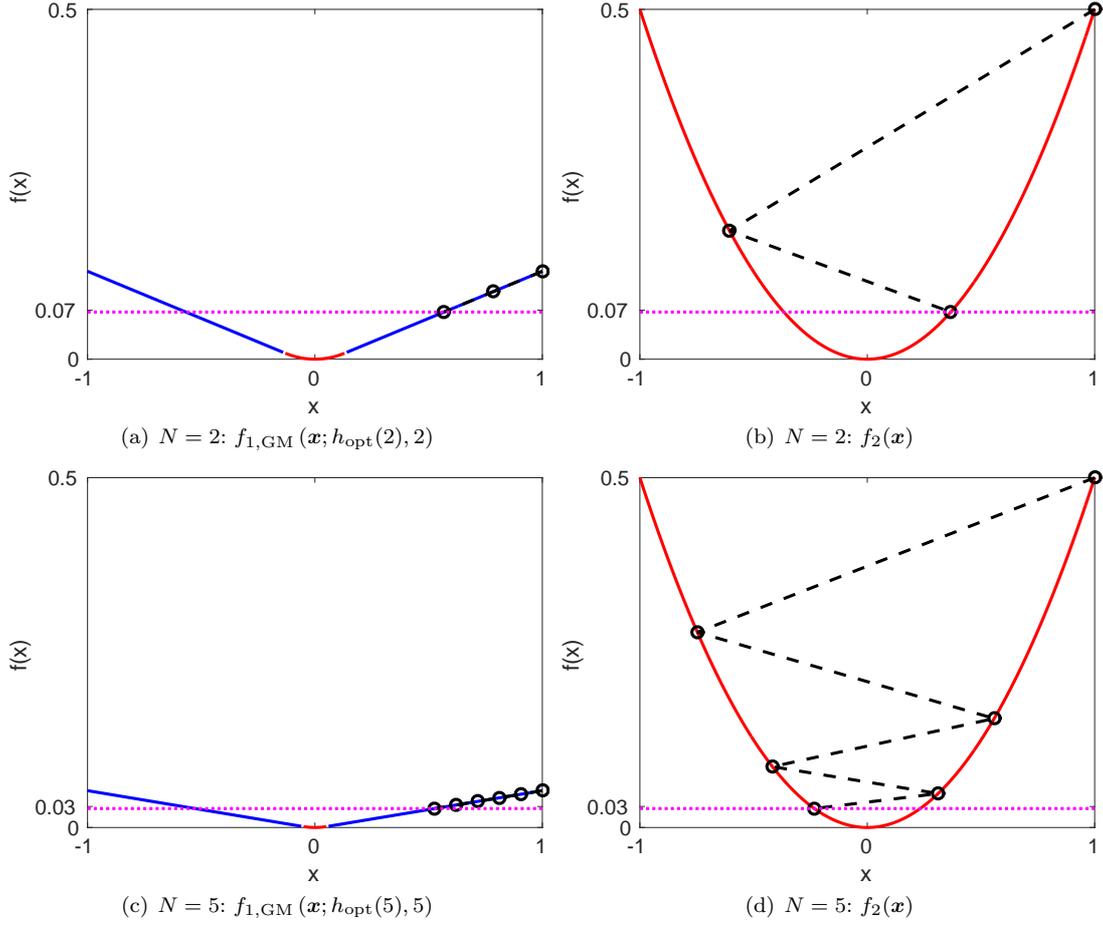

\begin{center}
	\subfigure[$N=2$: $\fwoneGMz\paren{\bmx;\hopt(2),2}$]
	{\includegraphics[clip,width=0.45\textwidth]{fig/fig,gm,constant,f1,N2.eps}}
        \subfigure[$N=2$: \fwtwo]
	{\includegraphics[clip,width=0.45\textwidth]{fig/fig,gm,constant,f2,N2.eps}}
	\subfigure[$N=5$: $\fwoneGMz\paren{\bmx;\hopt(5),5}$]
	{\includegraphics[clip,width=0.45\textwidth]{fig/fig,gm,constant,f1,N5.eps}}
        \subfigure[$N=5$: \fwtwo]
	{\includegraphics[clip,width=0.45\textwidth]{fig/fig,gm,constant,f2,N5.eps}}
\end{center}
\caption{The worst-case performance of the sequence $\{\bmx_i\}_{i=0}^N$ of GM 
	with an optimal 
	fixed-step \hoptN for $N=2,5$ and $d=L=R=1$.
The numerically optimized fixed-step sizes for $N=2,5$ 
are $\hopt(2) = 1.6058$ and $\hopt(5) = 1.7471$~\cite{taylor::ssc}.}
\label{fig:gm,constant}
\end{figure}

\subsection{Two worst-case functions for the last iterate $\bmx_N$ of OGM}
\label{sec:two,worst,ogm}

\cite[Theorem 3]{kim::ofo} showed that
\fwoneOGM~\eqref{eq:fwoneOGM} is a worst-case function 
for the last iterate $\bmx_N$ of OGM.
The following theorem shows that a quadratic function \fwtwo~\eqref{eq:fwtwo} 
is also a worst-case function for the last iterate of OGM.
\begin{theorem}
For the quadratic function $\fwtwo = \frac{L}{2}||\bmx||^2$~\eqref{eq:fwtwo} in \cF,
both OGM1 and OGM2 exactly achieve the convergence bound~\eqref{eq:ogm,conv}, \ie,
\begin{align*}
\fwtwoz(\bmx_N) - \fwtwoz(\bmx_*) = \frac{LR^2}{2\theta_N^2}
.\end{align*}
\begin{proof}
We use induction to show that
the following iterates:
\begin{align} 
\bmx_i = (-1)^i\frac{1}{\theta_i}R\bmnu, \quad i=0,\ldots,N,
\label{eq:quadpoint}
\end{align}
correspond to the iterates of OGM1 and OGM2 applied to \fwtwo.

Starting from $\bmx_0 = R\bmnu$,
where \bmnu is a unit vector,
and assuming that~\eqref{eq:quadpoint} holds for $i < N$,
we have
\begin{align*}
\bmx_{i+1} &= \bmx_i - \frac{1}{L}\sum_{k=0}^i\Hhkip\nabla\fwtwoz(\bmx_k) \\
	&= \paren{\bmx_i - \frac{1}{L}\Hhiip\nabla\fwtwoz(\bmx_i)} 
		- \frac{1}{L}\sum_{k=0}^{i-1}\frac{\theta_i - 1}{\theta_{i+1}}
			\hat{h}_{i,k}\nabla\fwtwoz(\bmx_k)
		+ \frac{1}{L}\frac{\theta_i-1}{\theta_{i+1}}\fwtwoz(\bmx_{i-1}) \\
	&= \frac{1 - 2\theta_i}{\theta_{i+1}}\bmx_i 
		+ \frac{\theta_i-1}{\theta_{i+1}}(\bmx_i - \bmx_{i-1})
		+ \frac{\theta_i-1}{\theta_{i+1}}\bmx_{i-1} 
	= - \frac{\theta_i}{\theta_{i+1}}\bmx_i \\
	&= (-1)^{i+1}\frac{1}{\theta_{i+1}}R\bmnu,
\end{align*}
where the second and third equalities use~\eqref{eq:fo} and~\eqref{eq:ogm1,h}.
Therefore, we have
\begin{align*}
\fwtwoz(\bmx_N) - \fwtwoz(\bmx_*) 
= \fwtwoz\paren{(-1)^N\frac{1}{\theta_N}R\bmnu} = \frac{LR^2}{2\theta_N^2}
\end{align*}
after $N$ iterations of OGM1 and OGM2, exactly matching the bound~\eqref{eq:ogm,conv}.
\qed
\end{proof}
\end{theorem}

Thus the last iterate $\bmx_N$ of OGM has two worst case functions:
\fwoneOGM and \fwtwo,
similar to an optimal fixed-step GM in Section~\ref{sec:two,worst,gm}.
Fig.~\ref{fig:worst,ogm} illustrates behavior of OGM for $N=2$ and $N=5$,
where OGM reaches same worst-case cost function value
for two different functions \fwoneOGM and \fwtwo
after $N$ iterations.

\begin{figure}[h!]
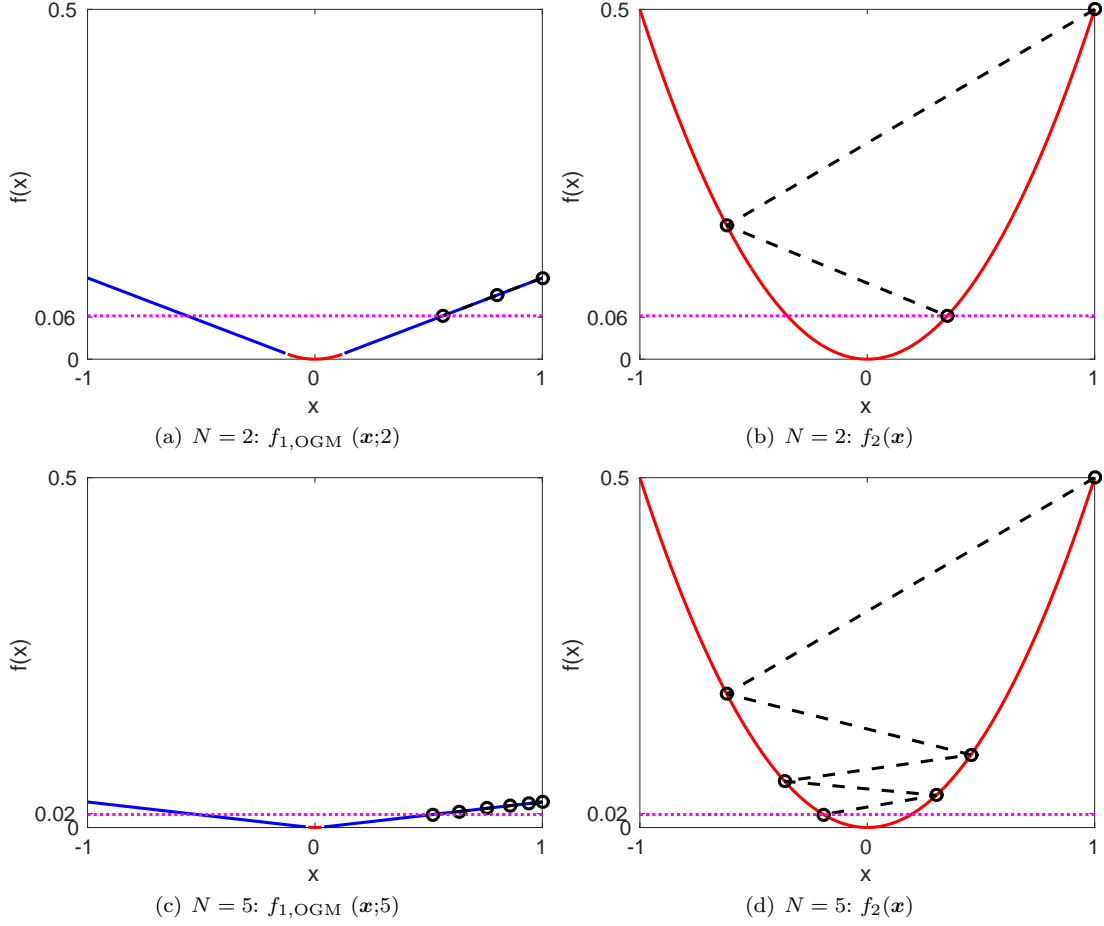

\begin{center}
        \subfigure[$N=2$: \fwoneOGMz(\bmx;2)]
        {\includegraphics[clip,width=0.45\textwidth]{fig/fig,ogm,f1,N2.eps}}
        \subfigure[$N=2$: \fwtwo]
        {\includegraphics[clip,width=0.45\textwidth]{fig/fig,ogm,f2,N2.eps}}
        \subfigure[$N=5$: \fwoneOGMz(\bmx;5)]
        {\includegraphics[clip,width=0.45\textwidth]{fig/fig,ogm,f1,N5.eps}}
        \subfigure[$N=5$: \fwtwo]
        {\includegraphics[clip,width=0.45\textwidth]{fig/fig,ogm,f2,N5.eps}}
\end{center}
\caption{The worst-case performance of the secondary sequence $\{\bmx_i\}_{i=0}^N$
        of OGM for $N=2,5$ and $d=L=R=1$.}
\label{fig:worst,ogm}
\end{figure}

In~\cite[Conjecture 4]{taylor::ssc} and Section~\ref{sec:conv,ogm,second},
the primary sequence of OGM 
is conjectured to have $\fwoneOGMp(\bmx;N)$ as a worst-case function,
whereas the quadratic function \fwtwo becomes the best-case
as the first primary iterate of OGM reaches the optimum just in one step.
On the other hand,
Section~\ref{sec:comp,worst} conjectured 
that \fwtwo is a worst-case function
for the secondary sequence of OGM
prior to the last iterate.
Apparently the primary and secondary sequences of OGM
have two extremely different
worst-case analyses,
whereas the last iterate $\bmx_N$ of OGM compromises between 
the two worst-case behaviors,
making the worst-case behavior of the optimal OGM interesting.

%

\section{Conclusion}
\label{sec:conc}

We provided an analytical convergence bound
for the primary sequence of OGM1 and OGM2,
augmenting the bounds of the last iterate 
of the secondary sequence of OGM in~\cite{kim::ofo}.
The corresponding convergence bound 
is twice as small as that of Nesterov's FGM,
showing that the primary sequence of OGM is 
faster than FGM.
However, interestingly the intermediate iterates of the \emph{secondary} sequence of OGM
were found to be slower than FGM in the worst-case.

We proposed two new formulations of OGM, called OGM1$'$ and OGM2$'$
that are related closely to
Nesterov's accelerated first-order methods in~\cite{nesterov:13:gmf}
(originally developed for nonsmooth composite convex functions
and differing from FGM in~\cite{nesterov:83:amf,nesterov:05:smo}).
For smooth problems, OGM and OGM$'$ provide faster convergence speed
than~\cite{nesterov:13:gmf}
considering the number of gradient computations required per iteration.

We showed that the last iterate of the secondary sequence of OGM 
has two types of worst-case functions,
a piecewise affine-quadratic function and a quadratic function.
In light of the optimality of OGM (for $d\ge N+1$) in~\cite{drori:16:tei},
it is interesting that
OGM has these two types of worst-case functions.
Because the optimal fixed-step GM
also appears to have two such
worst-case functions,
one might conjecture that this behavior is a general characteristic
of optimal fixed-step first-order methods.

In addition to the optimality of fixed-step first-order methods
for the cost function value,
studying the optimality for
an alternative criteria such as the gradient ($||\nabla f(\bmx_N)||$)
is an interesting research direction.
Just as Nesterov's FGM was extended for solving
nonsmooth composite convex functions~\cite{beck:09:afi,nesterov:13:gmf},
it would be interesting to extend OGM to such problems;
recently this was numerically studied by Taylor \etal~\cite{taylor:15:ewc}.
Incorporating a line-search scheme in~\cite{beck:09:afi,nesterov:13:gmf}
to OGM would be also worth investigating,
since computing the Lipschitz constant $L$ is sometimes expensive in practice.


\bibliographystyle{spmpsci}      
\bibliography{master0,mastersub}   

\end{document}